\documentclass[12pt]{article}
\title{Sur une question de Bergweiler
}

\font\sdopp=msbm10 scaled \magstep1
\def\CI {\sdopp {\hbox{C}}}
\author{Claudio Meneghini
\footnote{\tt 
Le courriel de l'auteur: clamengh@bluemail.ch
}
}
\font\cir=wncyb10

\def\pe{{\cir\hbox{P}}}

\begin{document}
\maketitle
\bibliographystyle{plain} 
\parindent=8pt
%caratteri cirillici
\font\cir=wncyb10
\def\Iu{\cir\hbox{YU}}
\def\Ze{\cir\hbox{Z}}
\def\pe{\cir\hbox{P}}
\def\Ef{\cir\hbox{F}}
\def\CIRC{\mathop{\tt o}\limits}
\def\quan{\vrule height6pt width6pt depth0pt}
\def\QUAN{\ \quan}
\def\BETA{\mathop{\beta}\limits}
\def\GAMMA{\mathop{\gamma}\limits}
\def\VI{\mathop{v}\limits}
\def\UI{\mathop{u}\limits}
\def\VII{\mathop{V}\limits}
\def\WI{\mathop{w}\limits}
\def\ZETA{\mathop{Z}\limits}

\newtheorem{definition}{Definition}%[chapter]
\newtheorem{defi}[definition]{D\'efinition}%[chapter]
\newtheorem{lemma}[definition]{Lemma}
\newtheorem{lemme}[definition]{Lemme}
\newtheorem{proposition}[definition]{Proposition}
\newtheorem{theorem}[definition]{Theorem}        
\newtheorem{theoreme}[definition]{Th\'eor\`eme}        
\newtheorem{corollary}[definition]{Corollary}  
\newtheorem{corollaire}[definition]{Corollaire}  
\newtheorem{remark}[definition]{Remark}  
\newtheorem{remarque}[definition]{Remarque}
  
%\font\sdopp=msbm10
\def\ERRE {\sdopp {\hbox{R}}}
\def\E {\sdopp {\hbox{E}}}
\def\QU {\sdopp {\hbox{Q}}}
\def\CI {\sdopp {\hbox{C}}}
\def\DI {\sdopp {\hbox{D}}}
\def\ENNE{\sdopp {\hbox{N}}}
\def\ZETA{\sdopp {\hbox{Z}}}
\def\ES{\sdopp {\hbox{S}}}
\def\PI {\sdopp {\hbox{P}}}
\def\M{\hbox{\tt\large M}}
\def\N{\hbox{\boldmath{}$N$\unboldmath}} 
\def\P{\hbox{\boldmath{}$P$\unboldmath}} 
\def\tr{\hbox{\boldmath{}$tr$\unboldmath}} 
\def\f{\hbox{\large\tt f}} 
\def\F{\hbox{\boldmath{}$F$\unboldmath}} 
\def\G{\hbox{\boldmath{}$G$\unboldmath}} 
\def\L{\hbox{\boldmath{}$L$\unboldmath}} 
\def\h{\hbox{\boldmath{}$h$\unboldmath}} 
\def\e{\hbox{\boldmath{}$e$\unboldmath}} 
\def\g{\hbox{\boldmath{}$g$\unboldmath}} 
\def\u{\hbox{\boldmath{}$u$\unboldmath}} 
\def\v{\hbox{\boldmath{}$v$\unboldmath}} 
\def\U{\hbox{\boldmath{}$U$\unboldmath}} 
\def\V{\hbox{\boldmath{}$V$\unboldmath}} 
\def\W{\hbox{\boldmath{}$W$\unboldmath}} 
\def\id{\hbox{\boldmath{}$id$\unboldmath}} 
\def\alph{\hbox{\boldmath{}$\alpha$\unboldmath}} 
\def\bet{\hbox{\boldmath{}$\beta$\unboldmath}} 
\def\gam{\hbox{\boldmath{}$\gamma$\unboldmath}} 
\def\pphi{\hbox{\boldmath{}$\varphi$\unboldmath}} 
\def\ppsi{\hbox{\boldmath{}$\psi$\unboldmath}} 
\def\Ppsi{\hbox{\boldmath{}$\Psi$\unboldmath}} 
\def\brevve{}
\def\Pert{}
\def\labelle #1{\label{#1}}
\def\quadras #1{ \hbox{\bf[}{#1}
\hbox{\bf]}}
\def\widepilde{}

\begin{abstract}
Nous montrons 
la densit\'e des cycles 
r\'epulsifs
dans 
l'ensemble de Julia des 
fonctions
m\'eromorphes
transcendantes 
\`a une variable complexe, sans 
utiliser le th\'eor\`eme
des cinq \^\i les d'Ahlfors
ni la th\'eorie de Nevanlinna.
\end{abstract}
%\Large
\section{Introduction}

Dans \cite{bergweiler} p.161,
Bergweiler posa la question s'il \'etait
possible de montrer la densit\'e des 
cycles r\'epulsifs dans l'ensemble
de Julia des 
fonctions
m\'eromorphes
transcendantes 
\`a une variable complexe, sans 
utiliser le th\'eor\`eme
des cinq \^\i les d'Ahlfors: 
Bolsch (voir \cite{bolsch}) a montr\'e que
la r\'eponse est oui, mais il fait appel
au th\'eor\`eme
des quatre valeurs compl\`etement ramifi\'ees.

On peut donc se poser la question si cette
d\'emonstration est encore possible
sans invoquer 
la th\'eorie de Nevanlinna

Le but de cette note-ci est de montrer que 
la r\'eponse est encore oui.
Tout d'abord, on
 utilisera, comme Bolsch, un th\'eor\`eme
 de Lehto sur la croissance
de la d\'eriv\'ee
sph\'erique d'une 
fonction m\'eromorphe
sur un disque \'epoint\'e,
ayant une singularit\'e essentielle
 isol\'ee au centre du disque
(voir \cite{lehto}
th.1).
L'utilisation de ce th\'eor\`eme-l\`a
sera tout-\`a-fait diff\'erente: en effet
cet \'enonc\'e sera combin\'e avec
des techniques de renormalisation \`a 
la Zalcman (voir \cite{zalcman})
et d'autre type,
\`a partir aussi d'un lemme m\'etrique de M.Gromov
(voir \cite{gromov}, p.256).

Si la fonction $ f  $
en consid\'eration a au moins deux poles
o\`u un pole qui n'est pas une valeur omise,
alors
la composition de $ f  $ \`a la source avec une famille
de 
contractions bien choisies nous permetra de
construire une application 
holomorphe limite sur le plan \'epoint\'e
et d'appliquer aux objets ainsi obtenus un
raisonnement semblable \`a 
\cite{bertelootduval} 
(voir aussi \cite{berteloot}, 
p.46)
pour les applications 
rationelles de $ \PI   $.
Cependant, on ne sera pas concern\'e avec
une famille non normale d'applications, mais on envisagera
la seule fonction $ f  $, au voisinage d'une 
singularit\'e essentielle
 isol\'ee.

Par contre, si $ f  $
n'a que un pole qui {\it est}
une valeur omise,
alors
la composition des it\'er\'ees $ f^{\circ n}  $
 \`a la source avec une bonne famille
de 
contractions nous permetra de
construire une application 
{\it enti\`ere} limite
et d'invoquer, dans ce nouveau contexte,
le raisonnement de
\cite{bertelootduval}, presque tel quel 
(voir aussi \cite{berteloot}
p.46).

\section{Pr\'eliminaires}

Rappelons ici quelques 
d\'efinitions classiques:
une fonction m\'eromorphe $ f  $ sur $ \CI  $
est dite {\it transcendante} si elle n'est 
pas une fraction rationnelle;
(voir par exemple \cite{bergweiler}):
{\it l'ensemble de Fatou 
$ {\cal F}_f  $}
de
$ f  $
est d\'efini comme 
l'ensemble des points au 
voisinage desquels
les it\'er\'ees de $ f  $ sont bien d\'efinies
et
forment une famille
normale de fonctions holomorphes;
{\it l'ensemble de Julia 
$ {\cal J}_f  $} est le 
compl\'ementaire de 
$ {\cal F}_f  $.

La condition que les it\'er\'ees 
$ \{ f^{\circ n}  \}  $
soyent bien
d\'efinies est g\'en\'eralement n\'ecessaire:
consid\'erons par exemple une fonction 
avec au moins
deux poles, ou bien un pole
n'\'etant pas une valeur omise;
l'orbite en arriere du point a l'infini est
non vide, ce qui entra\^\i ne, gr\^ace au th\'eor\`eme de Picard, qu'elle
est un ensemble infini.
Ainsi, 
l'ensemble des singularit\'es
essentielles de chaque $ f^{\circ n}  $
sera \'egalement infini.
Notons que cela comporte, gr\^ace au th\'eor\`eme de Montel,
la normalit\'e des 
it\'er\'ees sur les ouverts
o\`u elles sont d\'efinies.

Au contraire, si $ f  $ a un seul pole \'etant
une valeur omise, alors les it\'er\'ees sont 
toujours bien d\'efinies, mais la normalit\'e 
n'est pas garantie: dans ce cas, la d\'efinition
d'ensemble de Fatou et de Julia est semblable
a celle du cas des fractions rationnelles
(voir
par exemple
\cite{bergweiler}, p.153-155).

Rappelons maintenant la notion 
de {\it distance sph\'erique}
et la notion sous-jacente
de {\it d\'eriv\'ee sph\'erique}
d'une fonction m\'eromorphe.
La {distance sph\'erique}
$\sigma(z,z^{\prime})$
de deux points de la sph\`ere
de Riemann est d\'efinie comme
la distance euclid\'eenne
 de leurs projections
st\'er\'eographiques.

On voit ais\'ement que 
$\sigma(z,z^{\prime})=
{2|z-z^{\prime}|}/
{\sqrt{(1+|z|^2)
(1+|z^{\prime}|^2)}}
$ si $z,z^{\prime}\in\CI$, et
$\sigma(z,\infty)=
{2}/
{\sqrt{(1+|z|^2)}}
$.

\'Etant donn\'ee une fonction m\'eromorphe
au voisinage de $z\in\CI$,
on d\'efinit sa
{\it d\'eriv\'ee sph\'erique}
en $z$ comme
$$
g^{\sharp}(z):=\lim_{\zeta\to z}
 \frac{\sigma
(f(z),f(\zeta))}
{\sigma(z,\zeta)}.
$$
Cette limite existe et
on a $g^{\sharp}=(1/g)^{\sharp}$,
d'o\`u l'utilit\'e de cette notion
dans le domaine des
fonctions m\'eromorphes; 
par exemple on sait, gr\^ace
au th\'eor\`eme de Marty, que
la normalit\'e d'une famille
de fonction m\'eromorphes
est equivalente au fait que
la famille de
leurs d\'eriv\'ees sph\'eriques
 soit born\'ee.
 
Rappelons aussi les
\'enonc\'es du th\'eor\`eme
1 de \cite{lehto} et du 
lemme de Hurwitz (voir par exemple
\cite{berteloot}, p.8):
\begin{theoreme}
Soit $ v\in\CI  $, $ {\cal W}  $ un voisinage de
$ v   $ in $ \CI  $; soit $ g  $ une application
holomorphe (\`a valeurs en $ \PI   $) sur
$ {\cal W}\setminus\{ v  \}    $, ayant une singularit\'e essentielle \`a $ v  $ et 
$ g^{\sharp}  $ la d\'eriv\'ee sph\'erique de
$ g  $.
Alors
$ \limsup_{z\to v} \vert z-v     \vert\cdot g^{\sharp}(z)\geq 1/2 $.
\labelle{rapp}
\end{theoreme}
\begin{lemme}
\labelle{hurwitz}
Soit $ \Omega  $ une r\'egion
de $ \CI  $ et $ \{ h_n   \} $ une suite
d'applications holomorphes
d\'efinies sur $ \Omega  $ et
\`a valeurs en $ \PI   $: si 
$ h_n  $ converge uniform\'ement sur tout
compact de $ \Omega  $ vers une application 
holomorphe non constante
$ h:\Omega\to\PI   $ et 
celle-ci prend la valeur $ \alpha\in\PI   $, alors
$ h_n  $ prend aussi bien la valeur $ \alpha  $
pour tout $ n$ assez grand.
\end{lemme}
Rappelons enfin qu'un  {\it cycle}
d'une fonction $ f  $
est un point fixe de l'une de ses
it\'er\'ees $ f^{\circ n}  $, que son {\it multiplicateur}
est la d\'eriv\'ee de $ f^{\circ n}  $ a ce point
et que le cycle est {\it r\'epulsif} si son multiplicateur
est de module plus grand que $ 1  $.
\vskip.2cm
Le lemme suivant est connu comme
le {\tt lemme de l'espace 
m\'etrique}
(voir \cite{gromov}, p.256).

\begin{lemme} 
Soit $(X,d)$  un   
espace m\'etrique  complet et 
$M:X\rightarrow [0,+\infty)$ 
une fonction localement born\'ee.
Soit $\sigma>0$:  alors 
pour tout $u\in M^{-1}(0,+\infty)$
il
existe
$w\in X$
tel que:
{\tt (i)} $\displaystyle d(u,w)
         \leq 
\quadras{\sigma M(u)}^{-1}$;
{\tt (ii)}  $M(w)\geq M(u)$ et
{\tt (iii)} 
         $d(x,w)\leq 
\quadras{\sigma M(w)}^{-1}
\Rightarrow
M(x)\leq 2 M(w)$.
\labelle{metric}
\end{lemme}
{\sf D\'emonstration} 
Supposons par l'absurde
que le lemme soit faux: 
alors il existe 
$ u \in X $ tel que , pour tout
$w\in X$,
l'un au moins des \'enonc\'es {\tt (i)}, {\tt (ii)}
et
{\tt (iii)} est faux.
En particulier, $v_0:=u$ doit violer la condition {\tt (iii)}.
Donc on peut trouver $ v_1\in X  $
tel que $ M(v_1)>2 M(v_0)  $ mais 
$ d(v_1,v_0)\leq {1}/\sigma M({v_0})  $.
Cela
entra\^\i ne que {\tt (i)} et {\tt (ii)} 
sont vrais pour
$w=v_1$ et, par cons\'equent,
{\tt (iii)} doit \^etre faux pour $w=v_1$.
Donc on peut trouver $v_2\in X$ tel que 
$M(v_2)>2M(v_1)$
mais
$d(v_2,v_1)>[\sigma M(v_1)]^{-1}$, et
donc
$d(v_2,v_0)>\frac{1}{2}[\sigma M(v_0)]^{-1}$.
Ceci
entra\^\i ne que {\tt (i)} et {\tt (ii)} sont 
vrais pour
$w=v_2$ et, par cons\'equent, 
{\tt (iii)} doit \^etre faux m\^eme pour 
$w=v_2$.

En continuant ce proc\'ed\'e, on peut construire, par induction, une
 suite $\{v_n\}$ telle que $v_0=u$, 
$M(v_{n})>2 M(v_{n-1})>
2^{n}M(v_0)$ et 
$d(v_{n},v_{n-1})\leq
\displaystyle
{2^{1-n}
\left[
\sigma M(v_0))
\right]^{-1}
}$.
Cette suite-l\`a est de Cauchy:
en soit $ \lambda  $ la valeur limite.
On voit que $M$ n'est pas born\'ee
au voisinage de $\lambda   $:
c'est une contradiction.
\QUAN
\vskip.2cm
Le lemme suivant renormalise
(moyennant composition \`a la source
avec des contractions bien choisies)
 une famille d'applications
holomorphes (\`a valeurs dans la sph\`ere de Riemann).
La preuve ci-d\'ecrite 
peut se trouver en 
\cite{marseille}; 
voir aussi \cite{zalcman}.

\begin{lemme}
{\tt (Zalcman)}
Si une famille 
${\cal F}:=
\{f_{\alpha}\}$
d'applications holomorphes sur 
$\DI$, \`a valeurs en $ \PI  $,
 n'est pas normale 
sur aucun voisinage de
$v\in\DI$, alors
il existe des
suites 
$v_n\to v$,
%$\{r_n\}\subset\ERRE^+$, 
%avec 
$r_n\downarrow 0$,
$\{f_n\}\subset{\cal F}$
et une application holomorphe non constante
$h$ sur 
$\CI$, \`a valeurs en $ \PI  $,
telles que
$f_n(v_n+r_n z)$ est bien d\'efini
-pour $ n  $ assez grand- sur tout compact
de $ \CI  $ et
on y a, uniform\'ement, 
$f_n(v_n+r_n z)\to h$.%;
%en outre,
%la d\'eriv\'ee sph\'erique
%$h^{\sharp}$
%est born\'ee sur $\CI$.
%
\labelle{zalcman}
\end{lemme}
{\sf D\'emonstration} 
Gr\^ace \`a la non normalit\'e
\`a $v$,
 on peut trouver des suites
 $\xi_n\rightarrow v$ 
en $\DI$ et 
$\{f_n\}\subset{\cal F}$ 
telles que $
 f_n^{\sharp}(\xi_n)
\geq n^3$.
On peut supposer, sans nuire \`a la g\'en\'eralit\'e, que
$\{\xi_n\}$
soit contenu dans un  sous-ensemble 
ferm\'e
$X$ de $\DI$.

Pour tout $n$, appliquons le lemme
\ref{metric} \`a $X$
avec la m\'etrique 
euclid\'eenne, 
$M=f_n^{\sharp}$, 
$u=\xi_n$
et $\sigma=1/n$.
Appelons ($ \circ  $), ($ \circ  \circ    $) 
et ($ \circ  \circ  \circ    $) les cons\'equences des \'enonc\'es {\tt (i)}, respectivement
{\tt (ii)} {\tt (iii)} du lemme:
on obtient $v_n\in X$ tel que:
($ \circ  $) 
$
\displaystyle d(\xi_n,v_n)
          \leq 1/n^2$,
($ \circ  \circ    $) $f_n^{\sharp}(v_n)
\geq n^3$ et 
($ \circ  \circ  \circ    $) $\vert x-v_n     \vert
\leq 
{n}[{f_n^{\sharp}(v_n)}]^{-1}
\Rightarrow
f_n^{\sharp}(x)
\leq 2 f_n^{\sharp}
(v_n)$.

Posons maintenant
$r_n:=
[{f_n^{\sharp}(v_n)}]^{-1}
$ et
 $h_n(w):=f_n(v_n +r_n w)$.
Chaque $ h_n  $ est bien d\'efini sur
$ \DI(0,n)  $ car, gr\^ace \`a 
($ \circ $) et
($ \circ  \circ    $) ci-dessus,
$\,  v_n\to v\, $ $\, n r_n\leq 1/n^2  $.
La famille $\{h_n\}$ 
est normale,  car, gr\^ace \`a 
($ \circ  \circ  \circ    $)
$(h_n)^{\sharp}\leq 2$
 sur $\DI(0,n)$: gr\^ace au 
th\'eor\`eme d'Ascoli, on peut extraire de
$\{h_n\}$ une sous-suite
uniform\'ement 
convergente,
sur tout compact de
$ \CI  $, vers une fonction
m\'eromorphe
limite $h$
telle que
${h}^{\sharp}(0)
=\lim_{n\to\infty}
{h_n}^{\sharp}(0)=1$; cela 
prouve que
$h$ n'est pas 
constante.
%Finalement, par holomorphie,
%${h}^{\sharp}(z)
%=\lim_{n\to\infty}
%{h_n}^{\sharp}(z)\leq 2$ 
%pour tout $ z\in\CI  $.
\QUAN

\vskip.3cm
Le lemme suivant 'renormalise' une application
holomorphe (\`a valeurs dans la sph\`ere de Riemann) au voisinage d'une singularit\'e essentielle isol\'ee.
\begin{lemme}
Soient $ v\in\CI  $,
${\cal W}$ 
un voisinage ferm\'e
de $v$,
$ g  $ une application holomorphe,
\`a valeurs en $\PI$, 
sur 
${\cal W}\setminus \{v\}$, ayant une
singularit\'e essentielle 
%(isol\'ee)
 \`a $ v  $.
Alors
l'alternative suivante a lieu:
{\bf soit} il existe 
des
suites $v_n\to v$,
$\{ r_n  \}\subset\ERRE^+$, 
avec $r_n\to 0$,
telles que
$g (v_n+r_n z)$ est bien d\'efini sur
tout compact de $ \CI  $ 
et $g (v_n+r_n z)$
y converge
vers une application
holomorphe non constante
$h:\CI\to\PI$;
{\bf soit} 
il existe $ \zeta\in\CI  $ et
des
suites $v_n\to v$,
$\{ r_n  \}\subset\ERRE^+$, 
avec $r_n\to 0$,
telles que
$g (v_n+r_n z)$ est bien d\'efini sur 
tout compact de
$ \CI\setminus \{ \zeta  \}  $
et $g (v_n+r_n z)$
y converge
vers une application
holomorphe non constante
$h:\CI\setminus \{ \zeta  \}\to\PI    $.
%: dans ce
%dernier cas, $ h  $ a une singularit\'e essentielle
%\`a $ \zeta  $. 
%dont la d\'eriv\'ee sph\'erique
%$ h^{\sharp}  $
%est uniform\'ement born\'ee sur $ \DI  $.
%
\labelle{1zalcman}
\end{lemme}
{\sf D\'emonstration:} 
soit $\alpha_n   $ une suite de nombres r\'eels
positifs, avec $ \alpha_n\to 0  $.
Gr\^ace au 
th\'eor\`eme \ref{rapp},
on peut trouver une suite
 $\{\xi_n\}$ 
en ${\cal W}$ 
telle que 
$ \vert \xi_n -v     \vert=\alpha_n$ et
et
$\Pert  g^{\sharp}
(\xi_n)   \Pert
\geq [\alpha_n]^{-1}$
pour $ n  $ assez grand.
Posons
%$\varrho_n:=\vert \xi_n-v     \vert   $
%et
$X_n:={\cal W}\setminus
\DI(v, \alpha_n/4)
$;
pour tout $n,\in\ENNE$, appliquons 
le lemme
\ref{metric} \`a $X_n$ 
avec la m\'etrique 
euclid\'eenne, 
$M=
g^{\sharp} $, 
$u=\xi_n$
et $\sigma=8$.
Appelons ($ \diamondsuit  $), ($ \diamondsuit  \diamondsuit    $) 
et ($ \diamondsuit  \diamondsuit  \diamondsuit    $) les cons\'equences des \'enonc\'es {\tt (i)}, respectivement
{\tt (ii)} {\tt (iii)} du lemme:
cela fournit $v_{n}\in X_n$ 
tel que:
($\diamondsuit    $) 
$
\,\displaystyle \vert\xi_n-v_{n}\vert
          \leq \alpha_n/8\, $,
($ \diamondsuit  \diamondsuit    $) 
$\,\Pert g^{\sharp} 
(v_n)     \Pert
\geq [\alpha_n]^{-1}$ et \\ 
($ \diamondsuit  \diamondsuit  \diamondsuit    $) 
$
\ \left\{\vert x-v_n     \vert
\leq 
\displaystyle
{1}/{(8\Pert g^{\sharp} 
(v_n)\Pert)} \right\}
\Rightarrow
\Pert  g^{\sharp} 
(x)      \Pert
\leq 2\Pert  g^{\sharp} 
(v_{n})      \Pert
$.

Posons maintenant
$r_{n}:=\quadras{16g^{\sharp} (v_{n})}^{-1}$ et
$
h_{n}(z):=g(v_{n} +r_{n} z)$.
Gr\^ace \`a 
($\diamondsuit    $)  et
($\diamondsuit  \diamondsuit    $)
  ci-dessus, on a
$
\vert v_n-v     \vert
\leq \vert v_n-\xi_n\vert+
\vert \xi_n-v\vert
\leq 
\frac{9}{8}\alpha_n
$ et
$
\vert v_n-v     \vert
\geq \vert\xi_n-v\vert-
\vert \xi_n-v_n\vert
\geq 
\frac{7}{8}\alpha_n
$.
Ainsi,
$ 
|v_n+r_n z-v|\leq 
(9/8+|z|/16)\alpha_n\leq
2\alpha_n
  $ et
$ 
|v_n+r_n z-v|\geq 
(7/8-|z|/16)\alpha_n\geq
\frac{1}{4}\alpha_n
  $ ce qui entra\^\i ne que
$ v_n+r_n z\in X_n   $ et
que
$ h_{n}  $ est bien d\'efini sur
$ \DI$ pour $ n  $ assez grand.
La famille $\{h_n\}$ 
est normale sur $ \DI  $,  car, 
gr\^ace \`a 
($ \diamondsuit  \diamondsuit  \diamondsuit    $)  
$\Pert h_n^{\sharp} \Pert\leq 1/8$ 
sur $\DI$.
Gr\^ace au th\'eor\`eme d'Ascoli,
on peut en extraire 
une sous-suite
uniform\'ement 
convergente $ \{ h_{n_k}  \}  $,
sur tout compact de
$ \DI  $.
%, vers une limite $h$
%telle que
%$\Pert {h}^{\sharp}(0)\Pert
%=\lim_{n\to\infty}
%\Pert{h_n}^{\sharp}(0)\Pert
%=1/16$; cela 
%prouve que
%$h$ n'est pas 
%constante.
Or, \`a 
une nouvelle
 extraction pr\`es, 
l'alternative suivante a lieu:
soit $ (v-v_{n_k})/r_{n_k}\to\infty  $ soit 
$ (v-v_{n_k})/r_{n_k}  $ converge \`a un certain
nombre complexe,
que nous appellerons $ \zeta  $.

Notons que cela entra\^\i ne que,
$  h_{n_k}  $ est bien d\'efini sur tout compact
de $ \CI  $, resp. $ \CI\setminus   $

Une ult\'erieure alternative a lieu:
soit la famille $ \{ h_{n_k}  \}  $ est normale
sur $ \CI  $ (resp. $ \CI\setminus
\{ \zeta  \}  $), soit il existe 
$ \zeta^{\prime}\in\CI  $ (resp. 
$\zeta^{\prime}\in \CI\setminus
\{ \zeta  \}  $)
tel que $ \{ h_{n_k}  \}  $ n'est pas normale
au voisinage de $ \zeta^{\prime}$.
Dans le premier cas, la d\'emonstration
est termin\'ee moyennant extraction
d'une valeur limite $ h  $: cette function
jouit de la propri\'et\'e que 
$\Pert {h}^{\sharp}(0)\Pert
=\lim_{k\to\infty}
\Pert{h_{n_k}}^{\sharp}(0)\Pert
=1/16$; cela 
prouve que
$h$ n'est pas 
constante.

Dans le deuxi\`eme cas, le lemme \ref{zalcman} fournit
des suites $ \{ w_k  \}\in\CI  $ et $ \{ s_k  \}\in\ERRE^+  $
avec $ w_k\to\zeta^{\prime}  $ et
$ s_k\downarrow 0   $ telles que
$ h_{n_k}  (w_k+s_k z)  $ converge
(\`a extraction pr\`es)
 uniform\'ement
sur tout compact de 
$ \CI  $ (resp. $ \CI\setminus
\{ \zeta  \}  $) vers une fonction holomorphe 
non constante $ h: \CI\to \PI   $
(resp. $ h:\CI\setminus \{ \zeta  \}\to\PI   $).
Cela conclut la d\'emonstration, car
%d'un c\^ot\'e 
on a:
$ h_{n_k}  (w_k+s_k z)=g(v_{n_k}  +r_{n_k}  
w_k+r_{n_k}  s_k z)  $, $ v_{n_k}  +r_{n_k}  
w_k\to v  $ et $r_{n_k}  s_k \to 0   $; 
%de l'autre
%c\^ot\'e,  
%$\Pert {h}^{\sharp}(0)\Pert
%=\lim_{k\to\infty}
%\Pert{h_{n_k}}^{\sharp}(0)\Pert
%=1/16$; cela 
%prouve que
%$h$ n'est pas 
%constante.
 %par holomorphie,
%$\Pert {h}^{\sharp}(z)
%\Pert =\lim_{n\to\infty}
%\Pert {h_n}^{\sharp}(z)
%\Pert\leq 1/8$ pour tout 
%$ z\in\DI  $.
\QUAN

\section{Le r\'esultat principal
}
Soient maintenant:
$ f  $ une fonction m\'eromorphe non constante
sur 
$ \CI  $, ayant au moins deux poles, ou bien
un pole qui n'est pas une valeur
omise;
$ {\cal F}_f  $ et $ {\cal J}_f  $ les ensembles
de Fatou e Julia respectivement;
%(voir
%par exemple
%\cite{bergweiler}, p.153-155
 %pour la d\'efinition);
${\cal C}^+_f$ l'ensemble post-critique
de $f$ et ${\cal E}_f$ son
ensemble exceptionel (voir encore
\cite{bergweiler}, p.156).
Rappelons
que ${\cal C}^+_f$  est d\'efini
comme l'orbite de l'ensemble
critique de $f$: c'est un ensemble
d\'enombrable.
En outre, ${\cal E}_f$ est 
l'ensemble des points 
$\zeta\in\CI$
tels q	ue 
$\bigcup_{n=1}^{\infty} f^{-n}(\zeta)$
est un ensemble fini.
Gr\^ace aux th\'eor\`emes
de Montel et de Picard, cet ensemble-ci
peut contenir au plus deux points. 

Par ailleurs, 
on n'a pas forc\'ement
${\cal E}_f\subset 
{\cal F}_f $ pour les fonctions
transcendantes. 

\begin{theoreme}
Les cycles r\'epulsifs de $ f  $
sont denses dans 
$ {\cal J}_f  $.
\labelle{principal}
\end{theoreme}
{\sf D\'emonstration:} rappelons
que 
$\displaystyle
 {\cal J}_f =
\overline{O^-(\infty)}=
\overline{\bigcup_{n=0}^{\infty}f^{-n}(\infty)}
 $.
C'est un ensemble parfait (voir \cite{bergweiler}
p.154 et p.161).
Comme ${\cal C}^+_f\cup {\cal E}_f$ est
d\'enombrable, il suffit d'approcher
tout point 
$  p\in {\cal J}_f \setminus({\cal C}^+_f
\cup {\cal E}_f
)    \cap f^{-\lambda}(\infty)\quad  \lambda\in\ENNE  $.
\begin{lemme}
$ \bigcup_{l=0}^{\lambda }f^{-l}(\infty)  $
ne peut pas s'accumuler sur
$ p_{  }  $.
\labelle{aqui}
\end{lemme}
{\sf D\'emonstration:} supposons par l'absurde qu'il existe
une suite $ \{ p_{  \nu}  \} \subset
\bigcup_{l=0}^{  }f^{-l}(\infty)$
telle que $ p_{  \nu}\to p_{  }  $;
on peut en extraire une sous-suite 
$ q_{  \nu}\to p_{  }  $ telle que
$ \{ q_{  \nu}  \} \subset
f^{-n}(\infty)$, pour un $ 1\leq n\leq\lambda     $.
Alors:
{\tt\large A)} si  $ 1\leq n\leq \lambda -1   $,
$ f^{\circ n}  $ est holomorphe
\`a $ p_{  }  $, $ f^{\circ n}
(p_{  }  )\in\CI  $ mais $ f(q_{  \nu
})\equiv\infty  $ pour tout $ \nu\in\ENNE  $: c'est une contradiction; 
{\tt\large B)} si  $n=  \lambda   $,
$ f^{\circ n}  $ a un pole
\`a $ p_{  }  $, $ f^{\circ n}
(p_{  }  )=f(q_{  \nu
})\equiv\infty  $ pour tout $ \nu\in\ENNE  $: 
ceci entra\^\i ne $ f^{\circ n}\equiv\infty  $,
une contradiction.
\QUAN{\scriptsize\  (lemme \ref{aqui})}
\vskip.2cm
{\sf Fin de le d\'emonstration du th\'eor\`eme
\ref{principal}:}
Donc, gr\^ace au lemme \ref{aqui},
$ f^{\circ  \lambda +1 }  $ a une singularit\'e essentielle isol\'ee \`a $ p_{  }  $.

On peut alors appliquer le lemme \ref{1zalcman}
avec $ g:=f^{\circ  \lambda+1 }  $, $ v=p_{  }  $
et trouver un point $ \zeta\in\CI  $ et
des suites $  p_{ n}
\to p_{  }$,
 $ 
r_{ n}\downarrow 0
$ tels que
$f^{\circ \lambda+1}(p_{ n}+
r_{ n}z)$ converge 
uniform\'ement
sur tout compact (pour $ n\to\infty  $)
de $ \CI\setminus\{ \zeta  \}  $
vers une application
holomorphe non constante
$h_{  }: \CI\setminus\{ \zeta  \}  \to\PI$: 
%telle que, pour tout $z\in\DI$,
%$ h_{  }^{\sharp} (z) \leq 2$.
Donc $ h  $ est soit une application
enti\`ere $ h:\CI\to\PI  $, soit 
elle a une singularit\'e essentielle
\`a $ \zeta  $.

Gr\^ace au th\'eor\`eme de Picard,
$ h(\CI\setminus\{ \zeta  \})  $
rencontre
$ {\cal J}_f  $.

Soit donc $ U\subset\CI\setminus\{ \zeta  \}  $
 un ouvert tel que
$h(U)\cap {\cal J}_f \not=\emptyset  $: alors, 
on a, gr\^ace au th\'eor\`eme de Montel, 
 $\cup_{q\geq 1}
f^{\circ q}(h(U))\supset {\cal J}_f \setminus
{\cal E}_f$; ainsi,
 il existe
$ z_0\in U  $ et 
$ {\eta}\in\ENNE  $ tels que
$ p=f^{\circ {\eta}}\circ\widepilde h(z_0)  $; on peut supposer, sans nuire \`a la g\'en\'eralit\'e,
$\widepilde h(U)\subset\CI   $ et 
$ \widepilde h^{\prime}\not=0  $
sur $ U  $.

Or, 
$f^{\circ {\eta}}\circ
 f^{\circ \lambda }(p_{ n}+
r_{ n}z)
-(p_{ n}+
r_{ n}z)$
converge, 
 apres eventuelle extraction, 
vers
$f^{\circ {\eta}}\circ\widepilde h -p$, donc, 
le lemme de Hurwitz 
(lemme \ref{hurwitz})
%(voir par exemple
%\cite{berteloot}, p.8)
nous passe une suite de points
$\{ z_{ n}  \}\to z_0$
telle que 
$f^{\circ {\eta}}\circ
 f^{\circ \lambda }(p_{ n}+
r_{ n}z_{ n} )
=(p_{ n}+
r_{ n}z_{ n} )$:
ainsi les points 
$q_{n}:=p_{n}
+
r_{ n}z_{ n} 
$ forment
une suite $ q_{  n}\to v  $
de points p\'eriodiques 
de $ f  $.
Ces points sont r\'epulsifs (pour $ n  $ assez grand), puisque on a, d'un cot\'e,
$$
r_{ n}\cdot
\left(f^{\circ \eta+\lambda }    \right)^{\prime}
(p_{ n}+
r_{ n}z_{ n} )
\to
\left[\left(f^{\circ \eta}    \right)^{\prime}    
(\widepilde h(z_0))
\right]\cdot \widepilde h^{\prime}(z_0),
$$
et de l'autre cot\'e, $r_{ n}\to 0   $,
$ \widepilde h^{\prime}(z_0)\not=0  $ et 
$ \widepilde h(z_0)$ n'est  pas un point critique
de $ f^{\eta}  $: en effet, 
$ f^{\eta}\circ\widepilde  h(z_0)=p\not\in {\cal C}^+_f   $.
Cela conclut la d\'emonstration.
\QUAN

\vskip.3cm
Il reste \`a montrer le cas d'une fonction
$ f  $
m\'eromorphe ayant un pole \`a un point
qui est
une valeur omise par $ f  $: pour ce faire,
nous adaptons
la m\'ethode de renormalisations 
des it\'er\'ees d\'epeinte en
\cite{bertelootduval}.
Supposons, sans perte de g\'en\'eralit\'e,
que la singularit\'e soit plac\'ee \`a $ 0  $.
Alors les it\'er\'ees
$ \{ f^{\circ n}  \}  $
 sont bien d\'efinies partout
en $ \CI\setminus \{ 0  \}  $ et 
$ p\in{\cal J}_f  $ si et seulement si 
$ \{ f^{\circ n}  \}  $ n'est pas une famille normale
au voisinage de $ p  $.
On a encore:
\begin{theoreme}
Les cycles r\'epulsifs de $ f  $
sont denses dans 
$ {\cal J}_f  $.
\end{theoreme}
{\sf D\'emonstration:} 
appliquons le lemme \ref{zalcman}
\`a la famille $ \{ f^{\circ n}  \} $, 
avec $ v=p$:
cela fournit
des suites $ \{ p_{ n}  \}
\to p$ et
 $ \{ 
r_{ n}\}\downarrow 0
$ telles que
$\{f^{\circ n}(p_{ n}+
r_{ n}z)\}$ converge 
uniform\'ement
sur tout compact de $ \CI  $
vers une application
holomorphe non constante
$\widetilde h:\CI\to\PI$.
Comme ${\cal C}^+_f\cup {\cal E}_f$ est
d\'enombrable, 
on peut supposer, sans nuire \`a la g\'en\'eralit\'e,
$  p\in {\cal J}_f \setminus({\cal C}^+_f
\cup {\cal E}_f
)    $.

Soit maintenant
 $ U\subset\CI  $ un ouvert tel que 
$\widetilde h(U)\cap{\cal J}_f\not=0 $;
puisque $\cup_{q\geq 1}
f^{\circ q}(h(U))\supset {\cal J}_f \setminus
{\cal E}_f$,
 il existe
$ z_0\in U  $ et 
$ {\eta}\in\ENNE  $ tels que
$ p=f^{\circ {\eta}}\circ\widetilde h(z_0)  $; on peut supposer, sans nuire \`a la g\'en\'eralit\'e,
$\widetilde h(U)\subset\CI   $ et 
$ \widetilde h^{\prime}\not=0  $
sur $ U  $.

Or, 
$f^{\circ {\eta}}\circ
 f^{\circ  n  }(p_{ n}+
r_{ n}z)
-(p_{ n}+
r_{ n}z)$
converge, 
 apres eventuelle extraction, 
vers
$f^{\circ {\eta}}\circ\widetilde h -p$, donc, 
le lemme de Hurwitz 
nous passe une suite de points
$\{ z_{ n}  \}\to z_0$
telle que 
$f^{\circ {\eta}}\circ
 f^{\circ  n  }(p_{ n}+
r_{ n}z_{ n} )
=(p_{ n}+
r_{ n}z_{ n} )$:
ainsi les points 
$q_{n}:=p_{n}
+
r_{ n}z_{ n} 
$ forment
une suite $ q_{  n}\to v  $
de points p\'eriodiques 
de $ f  $.
Ces points sont r\'epulsifs 
comme au th\'eor\`eme \ref{principal}.
\QUAN

{\bf Remerciement:}
l'auteur tient \`a remercier le referee
pour beaucoup d'observations
et suggestions qui
ont permis d'am\'eliorer
la redaction de cet article.

\end{document}